\documentclass[oneside,reqno, twoside, 11pt]{amsart}
\usepackage[pdftex]{graphicx}
\usepackage{xcolor}
\usepackage[T1]{fontenc}
\usepackage{lmodern}
\usepackage[english]{babel}
\usepackage{mathrsfs}

\usepackage{upgreek}

\usepackage{amsfonts}
\usepackage{extpfeil}
\usepackage{verbatim}

\usepackage{BOONDOX-cal}

\usepackage{microtype}
\usepackage[new]{old-arrows}
\usepackage{amsmath}
\usepackage[all, cmtip]{xy}
\usepackage{tikz-cd}

\usepackage{amsthm}
\usepackage{mathtools}

\usepackage{bbm}
\usepackage{paralist}
\usepackage[T1]{fontenc}
\usepackage[colorlinks,citecolor=magenta,pagebackref=true,urlcolor=magenta, bookmarksdepth=3]{hyperref}
\usepackage[nodisplayskipstretch]{setspace}
\usepackage{breqn}

\usepackage{nicefrac}
\usepackage[centering,margin=2.5cm]{geometry}

\usepackage[font=small,labelfont=bf]{caption}
\usepackage[labelformat=simple]{subcaption}

\usepackage{epsfig}

\DeclareFontFamily{OT1}{pzc}{}
\DeclareFontShape{OT1}{pzc}{m}{it}{<-> s * [1.10] pzcmi7t}{}
\DeclareMathAlphabet{\mathpzc}{OT1}{pzc}{m}{it}

\makeatletter
\newtheorem*{rep@theorem}{\rep@title}
\newcommand{\newreptheorem}[2]{%
	\newenvironment{rep#1}[1]{%
		\def\rep@title{#2~\ref{##1}}%
		\begin{rep@theorem}}%
		{\end{rep@theorem}}}

\makeatother

\theoremstyle{plain}
\newtheorem*{thm*}{Theorem}
\newtheorem{thm}{Theorem}[section]
\newtheorem{mthm}{\bf Theorem}

\newreptheorem{mthm}{Main Theorem}
\newreptheorem{mcor}{Corollary}
\newreptheorem{mlem}{Lemma}
\newreptheorem{thm}{Theorem}

\newtheorem{cor}[thm]{Corollary}
\newtheorem{lem}[thm]{Lemma}

\newtheorem*{lem*}{Lemma}
\newtheorem{prp}[thm]{Proposition}
\newtheorem{conj}[thm]{Conjecture}

\theoremstyle{definition}
\newtheorem{dfn}[thm]{Definition}

\newtheorem{rem}[thm]{Remark}

\newcommand{\longtwoheadrightarrow}{\longrightarrow\hspace{-1.2em}\rightarrow\hspace{.2em}}

\newcommand{\ann}{\mathrm{ann}}

\renewcommand{\AA}{\mathcal{A}}
\newcommand{\BB}{\mathcal{B}}

\newcommand{\st}{\mr{st}}

\newcommand{\sbseq}{\subseteq}
\newcommand{\spseq}{\supseteq}

\newcommand{\n}{\frk{n}}

\newcommand{\vanish}[1]{}

\def\V{{\bf V}}

\def\ker{\mathrm{ker}}

\def\sbs\subset
\def\sbseq{\subseteq}

\def\({\left(}
\def\){\right)}
\def\no={\,{\,|\!\!\!\!\!=\,\,}}

\def\no={\,{\,|\!\!\!\!\!=\,\,}}
\def\sbseq{\subseteq}

\def\sbseq{\subseteq}
\def\sbs\subset
\def\spseq{\supseteq}

\newcommand{\xqedhere}[2]{%
	\rlap{\hbox to#1{\hfil\llap{\ensuremath{#2}}}}}

\newcommand{\cm}[1]{}

\newcommand\mbf[1]{\mathbf{#1}}

\newcommand\mr[1]{\mathrm{#1}}

\newcommand{\fld}{\mathbbm{k}}

\renewcommand\deg{\mr{deg}}

\newcommand\x{\mathbf{x}}

\DeclareMathOperator{\lk}{lk}

\newcommand{\term}{\textup{tm}}

\usepackage{xcolor}
\newcommand{\ord}{\textup{ord}}

\title[$p$-anisotropy on the moment curve for homology manifolds and cycles]{$p$-anisotropy on the moment curve for homology manifolds and cycles}
\author[Adiprasito]{Karim Alexander Adiprasito}
\address{Karim Adiprasito, Sorbonne Université and Université Paris Cité, CNRS, IMJ-PRG, F-75005 Paris, France}
\email{karim.adiprasito@imj-prg.fr}

\author[Hou]{Kaiying Hou}
\address{Kaiying Hou, Department of Mathematics, UC Berkeley, Berkeley, CA 94720, United States}
\email{kaiying@berkeley.edu}

\author[Kiyohara]{Daishi Kiyohara}
\address{Daishi Kiyohara, Department of Mathematics, 
Harvard University, 02138 Cambridge, United States}
\email{dkiyohara@math.harvard.edu}

\author[Koizumi]{Daniel Koizumi}
\address{Daniel Koizumi, Department of Mathematics, University of Texas at Austin, Austin, TX 78712, United States}
\email{daniel.koizumi@utexas.edu}

\author[Stephenson]{Monroe Stephenson}
\address{Monroe Stephenson, Technische Universität Berlin,
10623 Berlin, Germany}
\email{stephenson@campus.tu-berlin.de}

\date{\today}


\begin{document}
	
\begin{abstract}
We prove that the Gorensteinification of the face ring of a cycle is totally $p$-anisotropic in characteristic $p$.
In other words, given an appropriate Artinian reduction, it contains no nonzero $p$-isotropic elements. 
Moreover, we prove that the linear system of parameters can be chosen corresponding to a geometric realization with points on the moment curve.
In particular, this implies that the parameters do not have to be chosen very generically.
\end{abstract}
	
	\maketitle
	
	\newcommand{\AR}{\mathcal{A}}
	\newcommand{\BR}{\mathcal{B}}
		\newcommand{\IRL}{\mathcal{I}}
	\newcommand{\CR}{\mathcal{C}}
	\newcommand{\Mu}{M}
	\newcommand{\Soc}{\mathcal{S}\hspace{-1mm}\mathcal{o}\hspace{-1mm}\mathcal{c}}
	\newcommand{\Socl}{{\Soc^\circ}}
	\renewcommand{\n}{\mbf{n}}

\section{Introduction}
Biased pairing, and the more restrictive property of (total) anisotropy, is a replacement of the Hodge-Riemann relations beyond positivity in arbitrary characteristic. These relations were introduced by the first author in \cite{AHL} and Papadakis-Petrotou in \cite{PP}, respectively. It implies, when appropriately used, the Lefschetz property for face rings.

\begin{mthm}[\cite{AHL, APP, PP}]\label{mthm:gl}
Consider an arbitrary infinite field $\fld$, a simplicial cycle $\mu$ of dimension $d-1$ over $\fld$, and the associated graded commutative face ring~$\fld[\mu]$. 

Then there exists an Artinian reduction $\AR(\mu)$ and an element $\ell$ in $\AR^1(\mu)$ such that 
for every $k\le\nicefrac{d}{2}$, we have the \emph{hard Lefschetz property:} We have an isomorphism 
			\[\BR^k(\mu)\ \xrightarrow{\ \cdot \ell^{d-2k} \ }\ \BR^{d-k}(\mu). \]
\end{mthm}

Here, $\BR$ denotes the Gorensteinification of $\AR$. That is, $\BR$ is the quotient of $\AR$ by the annihilator of the fundamental class (also sometimes called the volume polynomial in this context). 

Essentially, biased pairings dictate that the Hodge-Riemann relations do not degenerate at certain subspaces; total anisotropy dictates that it does not degenerate at any of them. In \cite{APP}, it was proven that, given $\fld$ of characteristic two or $\fld=\mathbb{Q}$, we have total anisotropy with respect to a transcendental field extension. This proof relied on an earlier argument for spheres by Papadakis and Petrotou. In general characteristic, we only obtained biased pairing, specifically the nondegeneration at monomial ideals. 

Three key issues remained, put forward as questions in \cite{APP}:

\begin{enumerate}
\item Does anisotropy in characteristic 2 naturally extend to anisotropy in characteristic $p$?
\item Is it possible to restrict Artinian reductions for which anisotropy holds?
In particular, can we choose the linear systems of parameters so that they correspond to a geometric realization on the moment curve?
\end{enumerate} 

The main theorem in our paper answers a variation of the first question, and the second question affirmatively.
\begin{thm}[Total $pm$-anisotropy on the moment curve]\label{thm:ap}
Consider any field $\fld$ of characteristic $p$, any positive integer $m$, and any $(d-1)$-dimensional simplicial cycle $|\mu|$ over $\fld$.
Assume that the following condition is satisfied:
\begin{itemize}
    \item[$(*)$] For any face $\tau$ supported on $\mu$ such that $0<|\tau| < d$, the link satisfies $$\tilde{H}_{d-1-|\tau|}(\lk_\tau|\mu|,\fld) \cong \fld$$
    where $\tilde{H}_\bullet$ denotes reduced homology and $|\mu|$ denotes the support of the cycle $\mu$.
\end{itemize}
Then, for some field extension $\fld'$ of $\fld$, we have an Artinian reduction $\AR^*(|\mu|)$ that is $pm$-anisotropic, i.e., for every element $u\in \BR^k(\mu)$, $k\le \frac{d}{pm}$, we have 
\[u^{pm}\ \neq\ 0.\]
Moreover, the coefficients of the linear system of parameters can be chosen along the moment curve.
\end{thm}	

The additional condition on links can be dropped, but we feel it is less important to go into the technical details here; we refer to the paper \cite{APPHAL} for details. 

We remark that the first question had a similar partial answer in \cite{APPHAL} and independently \cite{KLS}, and our approach follows the methods laid out in the former, using a rather general observation on the inner workings of the anisotropy arguments of \cite{PP} and \cite{APP}. In \cite{KLS}, the answer is obtained using a more verbatim generalization of \cite{PP} to characteristic $p$, which, though it is possible to apply for purposes of this paper, is a bit too restrictive to do so easily.

\begin{rem}
    For any $\fld$-homology manifold $M$ of dimension $d-1$, $\lk_\tau M$ has the same homology as a $d-|\tau| - 1$-dimensional sphere \cite[p.2]{NovikSwartz}.
    Thus, Theorem \ref{thm:ap} applies for homology manifolds and, in particular, triangulated manifolds.
\end{rem}

For a simplicial homology sphere $\varDelta$ over $\mathbb{F}_2$, the anisotropy property of $\AA^*(\varDelta)$ with ground field $\mathbb{Q}$ was proved in \cite{KX}.
As an application of the main theorem, we generalize it to an arbitrary exponent.

\begin{thm}\label{thm:rat}
Let $\varDelta$ be a simplicial complex of dimension $d-1$, and let $m$ be a positive integer.
If $\varDelta$ is a homology sphere over $\mathbb{F}_p$ for all prime divisors $p$ of $m$, then for some field extension $\fld$ of $\mathbb{Q}$, we have an Artinian reduction $\AA^*(\varDelta)$ that is is $m$-anisotropic i.e. for any nonzero element $u\in \AA^k(\varDelta)$ of degree $k\le \frac{d}{m}$, we have 
\[u^m\neq0\in\AA^{mk}(\varDelta).\]
\end{thm}
\begin{cor}
Let $\varDelta$ be a simplicial homology sphere of dimension $d-1$ over $\mathbb{Z}$. Then for some field extension $\fld$ of $\mathbb{Q}$, we have an Artinian reduction $\AA^*(\varDelta)$ that is $m$-anisotropic for any $m\le d$.
\end{cor}

The paper is organized as follows.
Section \ref{sec:basics} introduces the relevant definitions.
Section \ref{sec:diff} proves a key lemma about the annihilator of certain monomials in the Gorensteinification $\BR^*(\mu)$.
Section \ref{sec:estab} establishes Theorem \ref{thm:ap} by considering certain derivatives of the degree.
Section \ref{sec:rat_co} proves Theorem \ref{thm:rat}.
Finally, Section \ref{sec:open_problems} discusses potential directions for future research.

\section{Preliminaries}\label{sec:basics}
	For the definitions of links, stars, face rings, and the Artinian reductions of the face rings, we refer the reader to  \cite{AHL} and \cite{AY}.
The following exposition on Gorensteinification and the degree map has appeared in \cite{APP}, but we include it here for the reader's convenience since these notions are extensively used in the rest of the paper.

\subsection{Gorensteinification}
Fix an infinite field $\fld$. 
Let $\varDelta$ be a simplicial complex.
Consider an Artinian reduction $\AR^\ast(\varDelta)$ of a face ring $\fld[\varDelta]$ with respect to a linear system of parameters $\Theta$. It is instructive to think of $\AR^\ast(\varDelta)$ as a geometric realization of $\varDelta$ in $\fld^d$, with the coefficients of $x_i$ in $\Theta$ giving the coordinates of the vertex $i$, recorded in a matrix $\V$.

Now, we pick a simplicial cycle 
\[\fld\mu\ \longhookrightarrow\ H_{d-1}(\varDelta;\fld)\]
and its dual map  $H^{d-1}(\varDelta;\fld)\twoheadrightarrow \fld\mu^\vee$. Via the canonical isomorphism
\[H^{d-1}(\varDelta;\fld)\ \cong\ \AR^d(\varDelta),\]
see \cite[Section 3.8]{AHL} and \cite{TW}, we have the quotient map $\AR^d(\varDelta) \twoheadrightarrow \fld \mu^\vee$. 
This quotient map then induces a pairing through
\[\AR^k(\varDelta)\ \times\ \AR^{d-k}(\varDelta)\ \to\ \AR^d(\varDelta) \to \fld \mu^\vee \cong \fld.
\]
$\BR^k_\mu(\varDelta)$ is then defined to be the quotient of $\AR^k(\varDelta)$ by the annihilator in the above pairing.
\begin{dfn}
    Let $L^k\subset \AR^k(\varDelta)$ denote the subspace of elements that have zero pairing with $\AR^{d-k}(\varDelta)$ and let $L = \bigoplus_k L^k$.
    The \textbf{Gorensteinification} of $\AR^*(\varDelta)$, denoted by $\BR^*_\mu(\varDelta)$, is defined to be
    $$
    \BR^*_\mu(\varDelta):=\AR^*(\varDelta)/L.
    $$
\end{dfn}
This does not depend on the simplicial complex $\varDelta$. 
A crucial property of $\BR^\ast_\mu(\varDelta)$ is that it is a Poincar\'e duality algebra \cite{APP}.
The following is immediate.
\begin{prp}
Consider a simplicial complex $\varDelta$ as above and a cycle $\mu$ in it. Then the restriction
\[\AR^\ast(\varDelta)\ \longtwoheadrightarrow\ \AR^\ast(|\mu|)\]
where $|\mu|$ denotes the support of $\mu$, that is, the minimal simplicial complex containing $\mu$, induces an isomorphism of Gorensteinifications. In particular, we have
\[\BR^\ast_\mu(\varDelta)\ \cong\ \BR^\ast_\mu(|\mu|).\]
\end{prp}
For ease of notation, we shall abbreviate $\BR^\ast(\mu):= \BR^\ast_\mu(|\mu|)$.
For simplicial complexes with unique homology cycles of dimension $(d-1)$, we set $\mu^\vee$ to be the unique cohomology cycle. In particular, this would be the case for orientable and connected manifolds.

\begin{rem}
    Performing Artinian reductions in finite fields can be problematic, so this subsection assumes an infinite field.
    In Theorem \ref{thm:ap}, $\fld$ is allowed to be a finite field while the Artinian reduction $\AR^*(|\mu|)$ and Gorensteinification $\BR^*(\mu)$ are done in an infinite field $\fld'$.
    A cycle over $\fld$ can also be considered as a cycle over $\fld'$ through the embedding $\fld\subset \fld'$, so the definitions in this subsection applies by working over $\fld'$.
\end{rem}

\subsection{Degree Map}
As explained in the last subsection, we have a surjection $\phi:\AR^d(\varDelta)\twoheadrightarrow \fld$.
Because $L^d = \ker\phi$, $\,\phi$ induces a map 
$$\deg: \BR^d(\mu) = \AR^d(\varDelta)/L^d \longrightarrow k$$
which we call the \textbf{degree map}.
Note that the degree map is readily described by the coefficients of the simplicial \emph{cycle}: it is enough to define it on cardinality $d$ faces $F$, as face rings are generated by squarefree monomials \cite[Section~4.3]{Lee} \footnote{Although Lee proves this only in characteristic zero, the argument goes through in general. We shall use his ideas several times in this paper for general characteristic, provided they apply. We note that all the ideas we use readily extend without any modification of the argument.}. For a cardinality $d$ face $F$, we have
\[\deg(\x_F)\ =\ \frac{\mu_F}{|\V_{|F}|},\]
where $\mu_F$ is the oriented coefficient of $\mu$ in $F$, and we fix an order on the vertices of $\mu$ and compute the sign with respect to the fundamental class, and the determinant $|\V_{|F}|$ of the minor $\V_{|F}$ of $\V$ corresponding to $F$.
For instance, if $\mu$ is the fundamental class of a manifold, we have canonically
\[\deg(\x_F)\ =\ \frac{\mr{sgn}(F)}{|\V_{|F}|}.\] 
In fact, it is possible give a degree formula not just for monomials $\x_F$'s that correspond to faces, but also for an arbitrary monomial $\x^\alpha$.
The following formula will be crucial for Section \ref{sec:estab}.
\begin{lem}[{\cite[Theorem~11]{Lee}}]\label{lee:formula}
We have
\[ \deg (\mathbf{x}^\alpha)\ =\ \sum_{F \text{ facet containing } \mathrm{supp}\ \alpha } \deg(\mathbf{x}_{F}) \prod_{i\in \mathrm{supp}\ F} [F-i]^{\alpha_i-1}\]
where $[F-i]$ is the volume element of $F-i$.
\end{lem}

The volume element $[F-i]$ is computed as follows: we first fix a general position vector $\rho$, then for all $F$ and $i\in F$, we replace the $i$-th column in the matrix $\V|_F$ with $\rho$ and set $[F-i]$ to be the determinant of the resulting matrix. 
It can be shown that the formula is independent from the choice of the general $\rho$ \cite[Section~4.5]{Lee}.

\section{Lemma on Annihilator}\label{sec:diff}
We prove an important lemma that allows us to write an element $u\in \BR^*(\mu)$ in a specific form so that the derivative calculation in Section \ref{sec:estab} can be done.
Fix a field $\fld$ and assume that Condition $(*)$ in Theorem \ref{thm:ap} applies to the cycle $\mu$. 
Let $\fld'$ be any field extension of $\fld$ that has an infinite number of elements.
The results presented in this section are for  $\fld'$, i.e., the Artinian reductions and Gorensteinifications considered are quotients of the face ring $\fld'[|\mu|]$.
\begin{lem}\label{anf}
    Let $\tau\in |\mu|$ be a $d-m-1$-dimensional face where $m$ is non-negative.
    Then the $m$-th graded piece $\ann^m_{\BB^*(\mu)}(x_\tau)$ of the annihilator is spanned by non-face monomials of $\st_\tau|\mu|$.
\end{lem}
\begin{proof}
    Let $I = I_{\st_\tau|\mu|}\subset \AR^*(|\mu|)$ be the non-face ideal of the star.
    We then have
    \begin{align*}
        &\AR^m(|\mu|)/I^m \cong \AR^m(\st_\tau|\mu|) \cong \AR^m(\lk_\tau|\mu|)\\
        \cong&\left(\tilde{H}_m(\lk_\tau|\mu|,\fld')\right)^\vee
        \cong \left(\tilde{H}_m(\lk_\tau|\mu|,\fld))\otimes_\fld \fld'\right)^\vee
    \end{align*}
    where the second isomorphism on the first line is obtained by repeatedly applying the cone lemma \cite[Lemma 3.2]{AHL},
    the first isomorphism on the second line is due to the Ishida complex that relates the graded pieces of the Artinian reduction with reduced cohomologies,
    and the second isomorphism on the second line is because $\fld'$ is flat over $\fld$.
    Thus, Condition $(*)$ ensures that $\dim_{\fld'} \AR^m(|\mu|)/I^m = 1$.
    Note that the quotient map $q: \AR^m(|\mu|) \to \BR^m(\mu)$ satisfies $q(I^m) \subset \ann^m_{\BB^*(\mu)}(x_\tau)$.
    We therefore obtain the following commutative diagram.
\[\begin{tikzcd}
	0&{I^m} & {\AA^m(\lvert\mu\rvert)} & {\AR^m(\lvert\mu\rvert)/I^m}&0 \\
	 0& {\ann^m_{\BB^*(\mu)}(x_\tau)} & {\BR^m(\mu)} & {\BR^m(\mu)/\ann^m_{\BB^*(\mu)}(x_\tau)} & 0
	\arrow["", from=1-1, to=1-2]
        \arrow["", from=1-2, to=1-3]
        \arrow["", from=1-3, to=1-4]
        \arrow["", from=1-4, to=1-5]
        \arrow["", from=2-1, to=2-2]
        \arrow["", from=2-2, to=2-3]
        \arrow["", from=2-3, to=2-4]
        \arrow["", from=2-4, to=2-5]
	\arrow["q"', from=1-2, to=2-2]
	\arrow["q"', from=1-3, to=2-3]
        \arrow["\tilde{q}"', from=1-4, to=2-4]
\end{tikzcd}\]
Note that both $\BR^m(\mu)/\ann^m_{\BB^*(\mu)}(x_\tau)$ and $\AR^m(\lvert\mu\rvert)/I^m$
are one dimensional, the former being so because $\BR^m(\mu)/\ann^m_{\BB^*(\mu)}(x_\tau)\cong \BB^d(\mu)$.
Because $q$ is a surjection, $\tilde{q}$ is too, so
$\tilde{q}$ is an isomorphism.
Applying the snake lemma immediately gives us that $q$ maps $I^m$ surjectively to the annihilator.
Because $I^m$ is spanned by non-face monomials of $\st_\tau|\mu|$, the annihilator must be too.
\end{proof}

\begin{prp}\label{eogt}
Let $u\in \BR^m(\mu)$ be a nonzero element of degree $m$.
There exists a face $\tau \in |\mu|$ of dimension $d-m-1$ such that $x_\tau\cdot u\neq0 \in \BR^d(\mu)$. For any such $\tau,$ there exists a face $\sigma \in lk_\tau{|\mu|}$ such that $u$ can be written as 
$$u = \lambda_\sigma x_\sigma + \sum_{\alpha} \lambda_\alpha x^\alpha $$
where $\lambda_\sigma \neq 0$ and $\text{supp}(x^\alpha) \not\in \text{st}_\tau |\mu|$ for all $\alpha$.

\end{prp}
\begin{proof}
The existence of such $\tau$ is gauranteed by the fact that $\BB^*(\mu)$ is a Poincar\'e duality algebra and that the face monomials span each graded piece of $\BB^*(\mu)$.
If we choose any complementary face $\sigma$ of $\tau$ such that the union $\sigma\cup \tau$ is a facet contained in the cycle $\mu$, we can see that $ux_\tau = \lambda_\sigma x_\tau x_\sigma$ holds for some scalar $\lambda_\sigma$ because $\BB^d(\mu)$ is one-dimensional.
This means that $ u - \lambda x_\sigma$ annihilates $x_\tau$ in $\BB^*(\mu)$.
The statement now follows from Lemma \ref{anf}.
\end{proof}

\section{Establishing full $p$-anisotropy}\label{sec:estab}
In this section we prove Theorem \ref{thm:ap}.
Let $\fld$ be any field of characteristic $p$ and let $\mu$ be a $d-1$-dimensional simplicial cycle that satisfies Condition $(*)$.
Let $v_1,...,v_n$ be the vertices of $|\mu|$.
We set $$\fld':=\fld(t_1,...,t_n),$$
which we will soon show is the field extension appearing in Theorem \ref{thm:ap}.
We define the linear system $\Theta\coloneqq\{\theta_1,\theta_2,...,\theta_d\}$
for the Artinian reduction
by
$$
\theta_i = t_1^ix_1 + t_2^ix_2 + \cdots+
t_n^ix_n.
$$
so that $\Theta$ corresponds to a geometric realization on the moment curve.
For the rest of this section, $\BB^*(\mu)$ refers to the Gorensteinification formed from the Artinian reduction
$$
\AA^*(|\mu|):=\fld'[|\mu|]/(\theta_1,...,\theta_d)
$$

Let us introduce some notations for convenience.
For a vertex $v$ in the support of $\mu$, we denote the associated partial differential by $\nicefrac{\mr{d}}{\mr{d}t_v}=:\partial_{v}$.
We write $f\approx g$ if there exists a nonzero function $h$ such that $f=h^p g$.
The notation is justified by the equivalence $\partial(f)\neq0\Longleftrightarrow\partial(g)\neq0$ with respect to any differential operator $\partial$, which follows from the equation $\partial(h^pg)=h^p\partial(g)$ in any field of characteristic $p$.

\begin{lem}[Full $p$-anisotropy]\label{lem:p}
For any nonzero element $u\in \BR^k(\mu)$, $k\le\frac{d}{p}$, we have 
\[u^p\neq0\in\BR^{pk}(\mu).\]
\end{lem}
\begin{proof}
Let us write $d=pk+\ell$ with $\ell\ge0$.
By the Poincar\'{e} duality of $\BR^\ast(\mu)$ there exists a squarefree monomial $x_\tau\in\BR^{pk+\ell-k}(\mu)$ such that $x_\tau u\neq0$.
For later use, we rewrite the face $\tau$ as the union of two disjoint faces $\xi$, $\iota$ with $|\xi|=pk-k$ and $|\iota|=\ell$.

We will prove that $u^p x_\iota\neq0\in\BR^d(\mu)$.
It is clearly sufficient to show $\partial_{\xi}(\deg(u^p x_\iota))\neq0$.
Using Proposition \ref{eogt}, we can express $u$ in the form
\[u=\lambda x_{\sigma}+\sum_{\alpha\not\in\st_{\tau}|\mu|}\lambda_\alpha x_{\alpha}\]
where $\sigma$ lies in $\lk_{\tau|\mu|}$.
Then we can write
\[\partial_{\xi}(\deg(u^p x_\iota))
    =\partial_{\xi}(\deg(\lambda^px_{\sigma}^px_\iota+\sum_{\alpha}\lambda_\alpha^p x_{\alpha}^px_\iota))
    =\lambda^p\partial_{\xi}(\deg(x_{\sigma}^px_\iota))+\sum_{\alpha}\lambda_\alpha^p\partial_{\xi}(\deg( x_{\alpha}^px_\iota))\]
where no mixed terms appear because of characteristic $p$.
Note that $\alpha\cup \iota\not\in\st_{\xi}|\mu|$ for each $\alpha$.
It follows that the differential $\partial_{\xi}$ acts trivially on all terms that arise in Lee's formula to compute $\deg(\x_\alpha^px_\iota)$.
Hence we have
\begin{align}
    \partial_{\xi}(\deg(u^p x_\iota))
    &=\lambda^p\partial_{\xi}(\deg(x_\sigma^px_\iota))
\end{align}
Therefore, we only have to show that \[\partial_{\xi}(\deg(x_\sigma^px_\iota))\neq0.\]
Moreover, since the faces $\xi$, $\sigma$, $\iota$ are complimentary and together make a unique facet $F\coloneqq \xi\cup \sigma\cup \iota$, the differential operator $\partial_\sigma$ acts trivially on all terms in Lee's formula applied to $\deg(x_\tau^px_\iota)$ except for the term corresponding to $F$, which we denote by $\term(F)$, so we have
\[\partial_{\xi}(\deg(x_\sigma^px_\iota))=\partial_\xi(\term(F)).\]
The term of our interest $\term(F)$ is given by
\[\term(F)=\frac{\mu_F}{\det(F)}\prod_{s\in F}[F-i]^{e(s)}, \,\,e(s)\coloneqq\begin{cases}-1&s\in\xi\\0&s\in\iota\\p-1&s\in\sigma\end{cases}.\]
In particular, we have $e(s)=0$ if $s\in\iota$ and $e(s)\equiv p-1\pmod p$ if $s\in\xi\cup \sigma$.
Then we find 
\begin{equation}\label{eq:triple}
    \term(F)\approx\frac{\mu_F}{\det(F)}\prod_{s\in\xi\cup\sigma}[F-s]^{p-1}.
\end{equation}
The task now has reduced to showing that the differential operator $\partial_\sigma$ acts on \ref{eq:triple} nontrivially.

The fact that we choose a geometric realization on the moment curve simplifies our calculation for \ref{eq:triple} and we find that
\begin{align*}
    \mu_F^{-1}\term(F)
    \approx\frac{1}{\det(F)}\prod_{s\in\xi\cup \sigma}[F-s]^{p-1}
    \approx& \frac{1}{\det(F)}\prod_{s\in\xi\cup \sigma}\left(\prod_{\substack{i,j\in \xi\cup\sigma\cup\iota \\ i\neq j\\ i,j\neq s}}(t_j-t_i)^{p-1}\right)\\
    \approx& \prod_{\substack{i,j\in \xi\cup\sigma\\i<j}}(t_j-t_i)^1\prod_{\substack{i\in\xi\cup\sigma\\ 
    j\in\iota}}(t_j-t_i)^0\prod_{\substack{i,j\in\iota\\ i<j}}(t_j-t_i)^{-1}\\
    \approx& \prod_{\substack{i,j\in \xi\cup\sigma\\i<j}}(t_j-t_i)\prod_{\substack{i,j\in\iota\\ i<j}}(t_j-t_i)^{-1}.
\end{align*}
Notice that the last term only depends on the variables which correspond to $\iota$.
Therefore,
\begin{align*}
    \mu_F^{-1}\partial_\sigma(\term(F))\approx \prod_{\substack{i,j\in\iota\\ i<j}}(t_j-t_i)^{-1}\partial_\xi\left(\prod_{\substack{i,j\in \xi\cup\sigma\\i<j}}(t_j-t_i)\right).
\end{align*}
Note that $\partial_\sigma$ does not act on $\mu_F$ because $\mu_F\in \fld\subset \fld(t_1,...,t_n)$.
The product that $\partial_\xi$ acts on is exactly the determinant of the Vandermonde matrix, so upon labeling vertices $\xi=[1,pk-k]$ and $\sigma=[pk-k+1,pk]$, it is given by
\[\prod_{1\le i<j\le pk}(t_j-t_i)=\sum_{\pi\in S_{pk}}\mr{sgn} (\pi)t_1^{\pi(1)}\cdots t_{pk}^{\pi(pk)}.\]
Then the differential operator $\partial_\xi$ acts nontrivially on every term that corresponds to a permutation $\pi\in S_{pk}$ such that $\pi(1),\cdots,\pi(pk-k)$ are not divisible by $p$.
This finishes the proof.
\end{proof}
The following corollary of Lemma \ref{lem:p} implies Theorem \ref{thm:ap}.
\begin{cor}
Let $\BR^\ast(\mu)$ have top degree $d$, and $u\neq0\in\BB^k(\mu)$.
Suppose that for a positive integer $m$ there exists some $n$ such that $\frac{d}{k}\ge pn\ge m$.
Then we have $u^m\neq0\in\BB^{mk}(\mu)$.
\end{cor}
\begin{proof}
Choose minimal $n$ satisfying the condition.
First, notice that it is enough to show $u^{pn}\neq0$ since $pn$ is at least $m$.
By $p$-anisotropy, it suffices to show $u^n\neq0$.
We can reduce to the case $m=1$ since $n<m$ would follow from $m>1$.
\end{proof}

\section{Application to rational coefficients}\label{sec:rat_co}
The goal of this section is to prove Theorem \ref{thm:rat}.
In order to do so, we establish an auxiliary result about the choice of a basis.
Our argument is similar to the one in \cite{KX}.
As in Section \ref{sec:estab}, 
we have a variable $t_i$ for each vertex $v_i$ in $|\mu|$
and the linear system $\theta_i$ is defined on the moment curve again.

For later use, we introduce the notation $\ord_p(f)$ for a rational polynomial $f$ and a prime $p$.
Any rational function $f\in\mathbb{Q}(t_1,\cdots,t_n)$ can be written in the form $f=a\cdot\frac{P(t_1,\cdots,t_n)}{Q(t_1,\cdots,t_n)}$ with $a\in\mathbb{Q}$ and $P,Q\in \mathbb{Z}[t_1,\cdots,t_n]$ where each of $P$ and $Q$ has integer coefficients whose gcd is one.
Then we write $\ord_p(f)\coloneqq |a|_p$ where $|a|_p$ is the $p$-adic valuation of $a$.

Let the subscript $\fld$ in the notation $\AR^*(|\mu|)_\fld$ and $\BB^{*}(\mu)_\fld$ indicate that the ground field used to construct the face ring is $\fld(t_1,\cdots,t_n)$.
Up to scalar multiplication any element $v\in\BB^{d}(\mu)_{\mathbb{Q}}$ can be written in the form
\[v=\sum_{|I|=d}\lambda_Ix_I\]
where $\ord_p(\lambda_I)\ge0$ for all $I$ and $\ord_p(\lambda_I)=0$ for some $I$. 
In this case we can associate to $v$ the element
\begin{equation*}
    [v]=\sum_{|I|=d}[\lambda_I]\cdot x_I
\end{equation*}
of $\BB^{d}(\mu)_{\mathbb{F}_p}$ where $[\lambda_I]$ denotes the residue class of $\lambda_I$.
Then Lee's formula implies the following relation about the degree maps on $\BB^{d}(\mu)_{\mathbb{F}_p}$ and $\BB^{d}(\mu)_{\mathbb{Q}}$:
\begin{equation}\label{eq:bracketdeg}
    [\deg(v)]=\deg([v]).
\end{equation}

We are ready to prove the following auxiliary lemma regarding the linear independence of a collection of monomials with respect to two different fields $\mathbb{F}_p$ and $\mathbb{Q}$.
\begin{lem}\label{lemma}
Let $\mu$ be a simplicial cycle.
Let $x_{\sigma_1},\cdots,x_{\sigma_s}$ be a collection of squarefree monomials of degree $k$ that are linearly independent in $\BB^k(\mu)_{\mathbb{F}_p}$.
Then the same collection of squarefree monomials $x_{\sigma_1},\cdots,x_{\sigma_s}$ are linearly independent viewed in $\BB^k(\mu)_{\mathbb{Q}}$.
\end{lem}
\begin{proof}
Suppose that there exist rational polynomials $\lambda_1,\cdots,\lambda_s$ such that $u\coloneqq\lambda_1x_{\sigma_1}+\cdots+\lambda_sx_{\sigma_s}=0$.
Here we may choose the coefficients so that $\ord_p(\lambda_i)\ge0$ for all $i$ and $\ord_p(\lambda_i)=0$ for at least one $i$.
Then the corresponding element $[u]=[\lambda_1]x_{\sigma_1}+\cdots+[\lambda_s]x_{\sigma_s}$ of $\BB^d(\mu)_{\mathbb{F}_p}$ satisfies $[u]x_\tau=0$ for all squarefree monomials $x_\tau$ of degree $d-k$ by equation \ref{eq:bracketdeg}.
By Poincar\'{e} duality, we conclude that $[u]=0$ so the same collection of monomials is linearly dependent in $\BB^d(\mu)_{\mathbb{F}_p}$.
\end{proof}

\begin{proof}[Proof of Theorem \ref{thm:rat}]
The field extension $\fld$ of $\mathbb{Q}$ appearing in Theorem \ref{thm:rat} refers to $\fld = \mathbb{Q}(t_1,...,t_n)$ and we will prove $m$-anisotropy in the ring $\BR^*(\varDelta)_{\mathbb{Q}}$ where $\varDelta$ denotes our homology sphere.
It suffices to show the statement when $m$ is a prime numbers $p$.
Because $\varDelta$ is a homology sphere over $\mathbb{F}_p$, it must also be a homology sphere over $\mathbb{F}_p(t_1,...,t_n)$.
By \cite[p.3]{NovikSwartz}, the Artinian reduction over $\mathbb{F}_p(t_1,...,t_n)$ is already Gorenstein, i.e.,
$$
\AR^*(\varDelta)_{\mathbb{F}_p} \cong \BR^*(\varDelta)_{\mathbb{F}_p}.
$$
Because being a homology sphere over $\mathbb{F}_p$ implies being a homology sphere over $\mathbb{Q}$ \cite[Lemma 2.1]{KX}, we likewise have  $\AR^*(\varDelta)_{\mathbb{Q}} \cong \BR^*(\varDelta)_{\mathbb{Q}}$.
By Lemma \ref{lemma}, there exists a common basis consisting of squarefree monomials $x_{\sigma_1},\cdots,x_{\sigma_s}$ for $\AA^k(\varDelta)_{\mathbb{Q}}$ and $\AA^k(\varDelta)_{\mathbb{F}_p}$ since these two vector spaces have the same dimension.
Any nonzero element $u\in\AA^k(\varDelta)_{\mathbb{Q}}$, up to scalar, admits a presentation of the form $u=\sum_{i=1}^s \lambda_ix_{\sigma_i}$ such that $\ord_p(\lambda_i)\ge0$ for all $i$ and $\ord_p(\lambda_i)=0$ for at least one $i$.
Then the corresponding element
\[[u]\coloneqq\sum_{i=1}^s [\lambda_i]x_{\sigma_i},\]
of $\AA^k(\varDelta)_{\mathbb{F}_p}$ is nonzero by the choice of basis.
If we write $d=pk+\ell$, there exists a face $\iota$ of dimension $\ell-1$ such that $x_\iota [u]^p\neq0$ by the proof of Lemma \ref{lem:p}.
By equation \ref{eq:bracketdeg} we conclude $x_\iota u^p\neq0$.
\end{proof}
\begin{rem}
    The reason that Theorem \ref{thm:rat} is only for homology spheres is because in the proof, we need $\AA^k(\varDelta)_{\mathbb{Q}}$ and $\AA^k(\varDelta)_{\mathbb{F}_p}$ to have the same dimension, which is not true for an arbitrary cycle due to potential torsion.
\end{rem}

\section{Open problems} \label{sec:open_problems}
The previous paper \cite{APP} asked for anisotropy using a specific linear system of parameters, as well as a generalization of anisotropy in characteristic 2 to $p$-anisotropy in characteristic $p$. Our paper concludes both of these questions. 
However, we pose the following open questions.

\begin{conj}[Total anisotropy on the moment curve in characteristic $p$]\label{cj:tot_ani}
Let $\fld$ be any field of characteristic $p$, $\mu$ any $(d-1)$-dimensional cycle over $\fld$, and the associated graded commutative face ring~$\fld[|\mu|]$. Then, for some transcendental field extension $\fld'$ of $\fld$, we have an Artinian reduction $\AR^*(|\mu|)$ that is anisotropic, i.e. for every nonzero element $u\in \BR^k(\mu)$, $k\le \frac{d}{2}$, we have 
\[u^2\ \neq\ 0.\]
	\end{conj}	
We attempted to tackle Conjecture \ref{cj:tot_ani} by the current approach of engineering clever linear systems and taking derivatives of the degree map. 
However, this approach failed and we were left wondering whether a different method was necessary or that Conjecture \ref{cj:tot_ani} was false all together.

We also pose the following conjecture.
Since there is some debate among the authors about which generality has a chance, we state the most pessimistic version: 
\begin{conj}
Consider a reduced Gorenstein standard graded algebra of Krull dimension, at least the socle dimension. Then a generic Artinian reduction over an appropriate field extension is anisotropic.
\end{conj}
\textbf{Acknowledgements.} We thank Stavros Papadakis, Vasso Petrotou and Geva Yashfe for helpful discussions.
The work was supported by the European Research Council under the European Union's Seventh Framework Programme ERC Grant agreement Horizon Europe ERC Grant number: 101045750 /
Project acronym: HodgeGeoComb, and was part of the REU program in Summer 2022. The fourth author acknowledges the support of the National Science Foundation in the form of a travel grant (NSF RTG grant \#1840190), which has enabled him to attend the REU program.


\end{document}